\documentclass{amsart}
\usepackage{amssymb,
amsmath,
amsthm,
graphicx,
amscd,
multind,
eufrak,
hyperref,
}
\theoremstyle{plain}
\newtheorem{thm}{Theorem}[section]
\newtheorem{prop}[thm]{Proposition}
\newtheorem{prb}[thm]{Problem}
\theoremstyle{definition}
\newtheorem{ex}[thm]{Example}
\newtheorem{re}[thm]{Remark}

\newcommand{\CC}{{\mathbb C}}
\newcommand{\RR}{{\mathbb R}}
\newcommand{\QQ}{{\mathbb Q}}
\newcommand{\ZZ}{{\mathbb Z}}
\newcommand{\NN}{{\mathbb N}}
{\begin{figure} \begin{center}}%
{\end{center} \end{figure}}

\newcommand{\diag}{\operatorname{diag}\nolimits}

\newcommand{\la}{\langle}
\newcommand{\ra}{\rangle}

\newcommand{\End}{\operatorname{End}\nolimits}
\newcommand{\GL}{\operatorname{GL}\nolimits}
\renewcommand{\Im}{\operatorname{Im}}
\newcommand{\SL}{\operatorname{SL}\nolimits}
\newcommand{\SO}{\operatorname{SO}\nolimits}

\newcommand{\liea}[1]{\mathfrak{#1}}

\newcommand{\tr}{\operatorname{tr}\nolimits}

\renewcommand{\Re}{\operatorname{Re}}
\newcommand{\Res}{\mathrm{Res}}
\newcommand{\Sp}{\mathrm{Sp}}
\renewcommand{\O}{\mathrm{O}}

\begin{document}

\title{Euclidean Distance degrees\\ of real algebraic groups}

\begin{abstract}
We study the problem of finding, in a real algebraic matrix group,
the matrix closest to a given data matrix. We do so from the
algebro-geometric perspective of Euclidean distance degrees.  We recover
several classical results; and among the new results that we prove is
a formula for the Euclidean distance degree of special linear groups.
\end{abstract}

\author[J.~Draisma]{Jan Draisma}
\address[Jan Draisma]{
Department of Mathematics and Computer Science\\
Technische Universiteit Eindhoven\\
P.O. Box 513, 5600 MB Eindhoven, The Netherlands;
and Centrum voor Wiskunde en Informatica, Amsterdam,
The Netherlands}
\thanks{JD is supported by a Vidi grant from
the Netherlands Organisation for Scientific Research (NWO)}
\email{j.draisma@tue.nl}

\author[Jasmijn Baaijens]{Jasmijn A. Baaijens}
\address[Jasmijn A. Baaijens]{
Department of Mathematics and Computer Science\\
Technische Universiteit Eindhoven\\
P.O. Box 513, 5600 MB Eindhoven, The Netherlands}
\email{j.a.baaijens@student.tue.nl}

\maketitle

\section{The distance to a matrix group}
\label{sec:distance}
Let $V$ be an $n$-dimensional real vector space equipped with a positive
definite inner product $(.|.)$, and write $\End(V)$ for the space
of linear maps $V \to V$. The inner product gives rise to a linear
map $\End(V) \to \End(V),\ a \mapsto a^t$ called transposition and
determined by the property that $(av|w)=(v|a^t w)$ for all $v,w \in V$,
and also to a positive definite inner product $\langle .,.  \rangle$
on $\End(V)$ defined by $\langle a,b \rangle:=\tr(a^t b)$. This inner
product enjoys properties such as $\langle a,bc \rangle=\langle b^t
a,c \rangle$. The associated norm $||.||$ on $\End(V)$ is called the
Frobenius norm. If we choose an orthonormal basis of $V$ and denote
the entries of the matrix of $a \in \End(V)$ relative to this basis
by $a_{ij}$, then $||a||^2=\sum_{ij} a_{ij}^2$. We will use the words
matrix and linear maps interchangeably, but we work without choosing
coordinates because it allows for a more elegant statement of some of
the results. For $a,b \in \End(V)$ and $u,v \in V$ we write $a \perp b$
for $\langle a,b \rangle=0$, and $v \perp w$ for $(v|w)=0$.

Let $G$ be a Zariski-closed subgroup of the real algebraic group $\GL(V)
\subseteq \End(V)$ of invertible linear maps. In other words, $G$ is a
subgroup of $\GL(V)$ characterised by polynomial equations in the matrix
entries. Then $G$ is a real algebraic group and in particular a smooth
manifold. The problem motivating this note is the following.

\begin{prb} \label{prob:closestmatrix}
Given a general $u \in \End(V)$, determine $x \in G$ that minimises the
squared-distance function $d_u(x):=||u-x||^2$.
\end{prb}

Here, and in the rest of this note, {\em general} means that whenever
convenient, we may assume that $u$ lies outside some proper,
Zariski-closed subset of $\End(V)$. Instances of this problem
appear naturally in applications. For instance, the {\em nearest
orthogonal matrix} plays a role in computer vision \cite{Horn86}, and
we revisit its solution in Section~\ref{sec:preserving}.  More or less
equivalent to this is the solution to the orthogonal Procrustes problem
\cite{Schonemann66}. For these and other matrix nearness problems we
refer to \cite{Higham89, Keller75}. More recent applications include
structured low-rank approximation, for which algebraic techniques are
developed in \cite{Ottaviani13}.

The bulk of this note is devoted to {\em counting} the number of critical
points on $G$ of the function $d_u$, in the general framework of the
{\em Euclidean distance degree} (ED degree) \cite{Draisma13c}. In
Section~\ref{sec:critical} we specialise this framework to matrix
groups. In Section~\ref{sec:preserving} we discuss matrix groups
preserving the inner product. In particular, we derive a conjecturally
sharp upper bound on the ED degree of a compact torus preserving the
inner product, revisit the classical cases of orthogonal and unitary
groups, and express the ED degree as the algebraic degree of a certain
matrix multiplication map. Then in Section~\ref{sec:notpreserving}
we discuss two classes of groups not preserving the inner product:
the special linear groups, consisting of all determinant-one matrices,
and the symplectic groups. For the former we determine the ED degree
explicitly. We conclude the note with a conjecture for the latter.

\section*{Acknowledgments}
We thank Pierre-Jean Spaenlehauer for his help with computing the ED
degree of $\Sp_6$, and Rob Eggermont, Emil Horobe\c t, and Hanspeter
Kraft for useful suggestions.

\section{The ED degree and critical equations}
\label{sec:critical}

As is common in the framework of ED degree computations, 
we aim to count the critical points of the function $d_u$
over the complex numbers, as follows. A point $x \in G$ is critical
for $d_u$ if and only if $(u-x) \perp a$ for all
$a$ in the tangent space $T_x
G \subseteq \End(V)$. As $G$ is an algebraic group, we have $T_x
G=x \cdot T_1 G=x \liea{g}$, where $T_1 G=\liea{g}$ is the tangent space of $G$
at the identity element $1$, i.e., the Lie algebra of $G$. So criticality
means that
\[ 0=\langle u-x,xb \rangle=\langle x^t(u-x),b \rangle \]
for all $b \in \liea{g}$. Hence, given $u$, we look for the solutions
of the {\em critical equations}
\begin{equation} \label{eq:crit}
x^t(u-x) \perp \liea{g} \text{ subject to } x \in G.
\end{equation}
The number of solutions $x \in G$ to \eqref{eq:crit} can vary with $u
\in \End(V)$. But if we set $V_\CC:=\CC \otimes_\RR V$, let $G_\CC
\subseteq \GL_\CC(V_\CC) \subseteq \End_\CC(V_\CC)$ be the set of
complex points of the algebraic group $G$, and extend $\la .,. \ra$ to
a symmetric $\CC$-linear form on $\End_\CC(V_\CC)$ (and {\em not} to a
Hermitian form!), then the number of solutions to \eqref{eq:crit} will
not depend on $u$, provided that $u$ is sufficiently general. Following
\cite{Draisma13c}, we call this number the {\em Euclidean distance degree}
(ED degree for short) of $G$. This number gives an algebraic measure
for the complexity of writing down the solution to the minimisation
problem~\ref{prob:closestmatrix}. We now distinguish two classes of
groups: those that preserve the inner product $(.|.)$ and those that
do not.

\section{Groups preserving the inner product}
\label{sec:preserving}

Assume that $(xv|xw)=(v|w)$ for all $x \in G$, so that $G$ is a subgroup
of the orthogonal group of $(.|.)$. Then all elements $x \in G$ satisfy
$x^t x=I$ and hence $||x||^2=n$, that is, $G$ is contained in the
sphere in $\End(V)$ of radius $\sqrt{n}$. As a consequence, $\liea{g}$
is contained in the tangent space at $1$ to that sphere, which equals
$1^\perp$. Hence the critical equations simplify to
\begin{equation} \label{eq:critcompact}
x^t u \perp \liea{g} \text{ subject to } x \in G.
\end{equation}
In other words, given the data matrix $u$ we seek to find all $x \in G$
that satisfy a system of linear homogeneous equations.  Alternatively,
we can write the critical equations as $u \in x \cdot \liea{g}^\perp$.
This proves the following proposition.

\begin{prop} \label{prop:degree} 
If $G$ preserves the inner product $(.|.)$, then the ED degree
of $G$ equals the degree of the multiplication map $G_\CC \times
\liea{g}_\CC^\perp \to \End_\CC(V_\CC),\ (x,s) \mapsto x \cdot s$. For
general real $u$, among the {\em real} pairs $(x,s)$ satisfying $xs=u$,
the one with the largest value of $\tr(s)$ is the one that minimises
$d_u(x)$.
\end{prop}

In other words, the ED degree counts the number of ways in which a
general matrix $u$ can be decomposed as a product of a matrix in $G_\CC$
and a matrix in $\liea{g}_\CC^\perp$. The last statement follows from
\[ d_u(x)=||u-x||^2=\tr(u^t u)-2 \tr(u^tx) + n=
\tr(u^t u) - 2 \tr(s) + n , \]
in which only the second term is not constant. 

\subsection*{Orthogonal groups}
If $G$ is the full orthogonal group of $(.|.)$, then $\liea{g}$
is the space of skew-symmetric matrices. Hence the decomposition of
Proposition~\ref{prop:degree} boils down to the classical {\em polar
decomposition}, where one writes a general matrix $u$ as $u=xs$ with $x$
orthogonal and $s$ symmetric. If $(x,s)$ is a solution, then
\[ s^2=s^t s=(x^{-1} u)^t (x^{-1} u)=u^t u, \]
a quadratic equation for $s$ that has $2^n$ real solutions for general
real $u$. Indeed, write
\[ u^t u=y \diag(\lambda_1,\ldots,\lambda_n) y^t \]
where $y$ is orthogonal and the $\lambda_i$ are the eigenvalues
of $u^t u$ (which are positive and distinct for general
$u$). Then any of the symmetric matrices
\[ s=y \diag(\pm \sqrt{\lambda_1},\ldots,\pm \sqrt{\lambda_n}) y^t \]
is a solution of the quadratic equation, and for each of these the matrix
$x=u s^{-1}$ is orthogonal, since
\[ x^t x = s^{-1} (u^t u) s^{-1} = 1. \]
We summarise our findings in the following, well-known theorem (see,
e.g., \cite{Keller75}).

\begin{thm} \label{thm:Orthogonal}
The ED degree of the orthogonal group of the $n$-dimensional inner
product space $V$ is $2^n$. Moreover, for general real $u$, all $2^n$
critical points of the squared distance function $d_u$ are real. The
critical point that minimises $d_u$ is
\[ x=u s^{-1} \text{ with }
 s=y \diag(\sqrt{\lambda_1},\ldots,\sqrt{\lambda_n})
y^t,
\]
where $y$ is an orthogonal matrix of eigenvectors of $u^tu$, and the
$\lambda_i$ are the corresponding eigenvalues.
\end{thm}

The last statement follows since $\tr(s)$ is maximised by taking the
{\em positive} square roots of the $\lambda_i$.

\begin{re}
The closest orthogonal matrix to a general real $u$ has determinant $1$
if $\det(u)>0$ and determinant $-1$ if $\det(u)<0$. Half of the critical
points of $d_u$ on the orthogonal group have determinant $1$, and half of
the points have determinant $-1$, that is, the ED degree of the special
orthogonal group is $2^{n-1}$. To find the special orthogonal
matrix closest to a matrix $u$ with $\det(u)<0$, one replaces the smallest $\sqrt{\lambda_i}$ in the construction above
by $-\sqrt{\lambda_i}$.
\end{re}

\subsection*{Unitary groups}

Assume that $n=2m$ and let $V$ be an $m$-dimensional complex vector
space, regarded as an $n$-dimensional real vector space. Let $h$
be a non-degenerate, positive definite, Hermitian form on $V$, where
we follow the convention that $h(cv,w)=ch(v,w)=h(v,\overline{c}w)$.
Define $(v|w):=\Re h(v,w)$. Then $(.|.)$ is a positive definite inner
product on $V$ regarded as a real vector space, and the norm on $V$
coming from $(.|.)$ is the same as that coming from $h$. Let $G$ be the
unitary group of $h$, which consists of all $x:V \to V$ that are not
only $\RR$-linear but in fact $\CC$-linear and that moreover preserve
$h$. Such maps $x$ also preserve $(.|.)$, so we are in the situation
of this section. The converse is also true: if $x$ is
$\CC$-linear and preserves $(.|.)$, then
\begin{align*}
\Im h(u,v)&=-\Re
h(iv,w)=-(iv|w)=-(x(iv)|x(w))=-(ix(v)|x(w))\\
&=-\Re h(ix(v),x(w))=\Im h(x(v),x(w)),
\end{align*}
so that $x$ preserves both $\Re h$ and $\Im h$ and hence $h$.

Note that $\End_\RR(V)$ has dimension $n^2=4m^2$, but $G$ is contained
in the real subspace $\End_\CC(V) \subseteq \End_\RR(V)$, which has
real dimension $2m^2$. For a general data matrix $u \in \End_\RR(V)$,
the critical points of $d_u$ on $G$ will be the same as the critical
points of $d_{u'}$ where $u'$ is the orthogonal projection of $u$ in
$\End_\CC(V)$. Hence in what follows we may assume that $u$ already
lies in $\End_\CC(V)$, and we focus our attention entirely on the space
$\End_\CC(V)$. For a linear map $u$ in the latter space, we write $u^*$
for the $\CC$-linear map determined by $h(v,uw)=h(u^*v,w)$ for all $v,w
\in V$. This map will also have the property that $(v|uw)=(u^*v|w)$,
i.e., $u^*$ {\em coincides} with our transpose $u^t$ relative to $(.|.)$.
In what follows we follow the convention to write $u^*$.

The Lie algebra $\liea{g}$ consists of all skew-Hermitian linear maps
in $\End_\CC(V)$, and its orthogonal complement $\liea{g}^\perp$ inside
$\End_\CC(V)$ therefore consists of all Hermitian $\CC$-linear maps.
Again, the decomposition of Proposition~\ref{prop:degree} boils down to
the polar decomposition. Here is the result (see, e.g., \cite{Keller75}).

\begin{thm}
The ED degree of the unitary group of a non-degenerate Hermitian form
$(.|.)$ on an $m$-dimensional complex vector space $V$ equals $2^m$. For
a general data point $u \in \End_\CC(V)$ the critical points
are computed as follows. First write
\[ u^* u=y \diag(\lambda_1,\ldots,\lambda_m) y^*, \]
where $y$ is a unitary map and the $\lambda_i \in \RR_{\geq 0}$ are the
eigenvalues of $u^* u$, then pick any of the $2^m$ square
roots
\[ s=y \diag(\pm \sqrt{\lambda_1},\ldots,\pm \sqrt{\lambda_m}) y^* \]
of $u^* u$, and finally set $x:=us^{-1}$. Choosing all square roots positive
leads to the closest unitary matrix to $u$.
\end{thm}

There is a slight subtlety in the last statement: to find the closest
matrix, we have to maximise $\tr_\RR(s)$, where we see $s$ as an
element of $\End_\RR(V)$, while the sum of the $m$ eigenvalues $\pm
\sqrt{\lambda_i}$ equals $\tr_\CC(s)$. But in fact, for any $z \in
\End_\CC(V)$, we have $\tr_\RR(z)=\tr_\CC(z) + \overline{\tr_\CC(z)}$,
so that $\tr_\RR(s)=2\tr_{\CC}(s)$.

\subsection*{Compact tori}

Assume that the real algebraic group $G$ is a compact torus, i.e.,
that it is an abelian compact Lie group and abstractly isomorphic to
a power $(S^1)^m$ of circle groups. We continue to assume that $G$
preserves the inner product $(.|.)$. We will bound the
ED degree of $G$, and unlike in the previous two examples we will make
extensive use of the complexification $G_\CC$ of $G$.

Indeed, $G_\CC$ is now isomorphic to an algebraic torus
$T:=(\CC^*)^m$. Let $X(T)$ be the set of all characters of $T$, i.e.,
of all algebraic group homomorphisms $\chi: T \to \CC^*$. These are all
of the form $\chi:(t_1,\ldots,t_m) \mapsto t_1^{a_1} \cdots t_m^{a_m}$,
where $a_1,\ldots,a_m \in \ZZ$; and this gives an isomorphism $X(T)
\cong \ZZ^m$ of finitely generated Abelian groups: $X(T)$ with respect
to multiplication and $\ZZ^m$ with respect to addition. We will identify
these groups, and accordingly write $t^\chi$ instead of $\chi(t)$ and
write $+$ for the operation in $X(T)$, so $t^{\chi+\lambda}=t^{\chi}
\cdot t^{\lambda}$ and $t^\chi=t_1^{a_1} \cdots t_m^{a_m}$
if $\chi=(a_1,\ldots,a_m)$.

The isomorphism $T \mapsto G_\CC \subseteq \GL_\CC(V_\CC)$ gives $V_\CC$
the structure of a $T$-representation. As such, it splits as a direct
sum of one-dimensional $T$-representations:
\[ V_\CC=\bigoplus_{\chi \in X(T)} V_\chi, \]
where, for any $\chi \in X(T)$, we let $V_\chi$ be the corresponding
eigenspace (or {\em weight space}), defined by
\[ V_{\chi}:=\{v \in V_\CC \mid \forall t \in T: tv=t^\chi v\}. \]
Of course, only finitely many of these spaces are non-zero, and their
dimensions add up to $n$. Let $X_V \subseteq X(T) = \ZZ^m$ denote the
set of characters $\chi$ for which $V_\chi$ is non-zero.  The fact that
the map $T \to G$ is an isomorphism implies that the lattice in $\ZZ^m$
generated by $X_V$ has full rank. We will prove the following result.

\begin{thm}
The ED degree of the compact torus $G \cong (S^1)^m$ depends only
on $X_V$, and is independent of the dimensions of the weight spaces
$\dim V_\chi,\ \chi \in X_V$. Moreover, it is bounded from above by the
normalised volume of the convex hull $\Delta$ of $X_V \subseteq \ZZ^m$.
Here the normalisation is such that the simplex spanned by $0$ and the
standard basis vectors has volume one.
\end{thm}

For the proof, observe that
the (complexified) bilinear form $(.|.)$ on $V$ is preserved by $T$. For
$v \in V_\chi$ and $w \in V_\lambda$ and $t \in T$ we
therefore have
\[ t^{\chi+\lambda} \cdot (v|w)=(t^\chi v|t^\lambda w)=
(t v | t w)=(v|w). \]
Hence, if $v$ and $w$ are not perpendicular, then $\chi + \lambda$ is
the trivial character sending all of $T$ to $1$ (so $\chi+\lambda=0 \in
\ZZ^m$). In other words, $(.|.)$ must pair each $V_\chi$ non-degenerately
with the corresponding space $V_{-\chi}$, and is zero on all pairs
$V_\chi \times V_\lambda$ with $\lambda \neq -\chi$. In particular,
we have $\dim V_\chi=\dim V_{-\chi}$ for all $\chi \in
X(T)$, and $X_V$ is centrally symmetric.

We can now choose a basis $v_1,\ldots,v_n$ of $V_\CC$
consisting of $T$-eigenvectors
such that $(v_i|v_j)=\delta_{j,n+1-i}$.  Let $\chi_i \in
X_V$ be the character of $v_i$, i.e., we have $t v_i=t^{\chi_i} v_i$
for all $t \in T$. There will be repetitions among the $\chi_i$ if some
of the weight spaces have dimensions greater than $1$; and by the above we
have $\chi_i=-\chi_{n+1-i}$. Relative to this basis we have
\[ G_\CC=\{(t^{\chi_1},\ldots,t^{\chi_n}) \mid t \in (\CC^*)^m\}
\text{ and }
\liea{g}_\CC=\{\diag(w \cdot \chi_1,\ldots, w\cdot \chi_n) \mid
w \in \CC^m\} \]
where $\chi_i \cdot w$ is the ordinary dot product of $w \in \CC^m$
with $\chi_i \in \ZZ^m$. The bilinear form $\langle .,.\rangle$ on
matrices takes the form
\[ \langle E_{ij}, E_{kl} \rangle=\delta_{k,n+1-i}
\delta_{l,n+1-j}, \]
where $E_{ij}$ is the map whose matrix relative to the basis
$v_1,\ldots,v_n$ has a $1$ at position $(i,j)$ and zeroes
elsewhere. Hence, for a data matrix $u \in \End(V)$, the critical
equations \eqref{eq:critcompact} translate into the following equations
for the pre-image $t \in T$ of $x \in G_\CC$:
\[ (t^{-\chi_1} u_{11})(w \cdot \chi_n) + \ldots +
(t^{-\chi_n} u_{nn}) (w \cdot \chi_1) =0 \text{ for all } w \in \CC^m.
\]
Using $\chi_i=-\chi_{n+1-i}$ we may rewrite this as
\begin{equation} \label{eq:BKK}
\sum_{i=1}^n (t^{\chi_i} u_{ii}) (w \cdot \chi_i) =0
\text{ for all } w \in \CC^m \end{equation}
The ED degree of $G$ is the number of solutions $t \in (\CC^*)^m$ to
this system of $m$ Laurent-polynomial equations for general values of the
$u_{ii}$. Grouping the indices $i$ for which the $\chi_i$ are equal and
adding up the corresponding $u_{ii}$ we find that the cardinality of the
solution set is, indeed, independent of the dimensions $\dim V_{\chi_i}$,
the first statement in the theorem.

Letting $w$ run over a basis of $\CC^m$, we obtain a system of $m$ Laurent
polynomial equations for $t \in (\CC^*)^m$ with fixed support $X_V$. The
Bernstein-Kushnirenko-Khovanskii theorem \cite[Theorems A,B]{Bernstein76}
ensures that the number of isolated solutions to this system is less
than or equal to the normalised volume of $\Delta$. This proves that
the ED degree does not exceed that bound, and hence the second part of
the theorem.

\begin{ex}
Let the group $\SO_2 \cong S^1$ act on $V:=(\RR^2)^{\otimes d}$
via $g(v_1 \otimes \cdots \otimes v_d)=(gv_1) \otimes \cdots \otimes
(gv_d)$, and let $(.|.)$ be the inner product on $V$ induced from that on
$\RR^2$. Let $G \subseteq \GL(V)$ be the image of $\SO_2$. Then the ED
degree of $G$ equals $2d$, computed as follows. All elements of $\SO_2$
are diagonalised over $\CC$ by the choice of basis $f_1:=e_1 + i e_2$
and $f_{-1}:=e_1 - i e_2$ of $\CC^2$. The complexification of $\SO_2$
is the image of the one-dimensional torus $T=\CC^*$ in its action on
$\CC^2$ via $(t,f_i) \mapsto t^i f_i$. The complexification $G_\CC$
is the image of $T$ in its induced action on $\CC^{\otimes d}$ via
\[ (t,f_{i_1} \otimes \cdots \otimes f_{i_d}) \mapsto
	t^{\sum_j i_j} f_{i_1} \otimes \cdots \otimes
	f_{i_d} \text{ for all } i_1,\ldots,i_d \in \{\pm
	1\}.  \]
We have $X_V=\{-d,-d+2,\ldots,d\} \subseteq \ZZ^1$, where $\ZZ^1$ is the
character lattice of $T$. If $d$ is odd, then the map $T \to G_\CC$
is one-to-one, and $\Delta$ is a line segment of length $2d$. If $d$
is even, then that map is two-to-one (since all exponents $\sum_j i_j$
above are then even), and the character lattice of $G_\CC$ is then
$2\ZZ^1 \subseteq \ZZ^1$. In this case, the normalised volume of $\Delta$
is $d$. The theorem says that the ED degree is at most $2d$ for odd $d$
and at most $d$ for even $d$. In this case equality holds: the system
\eqref{eq:BKK} reduces to
\[ u'_{-d} t^{-d} + \cdots + u'_{d} t^d = 0, \]
where the $u'_j$ are sums of $u_{ii}$ corresponding to the same character;
and for general values of the $u'_j$ this equation has exactly $2d$
solutions for $t$, and exactly $d$ solutions for $t^2$ when $d$ is even.
\end{ex}

We do not know if the system \eqref{eq:BKK}, for general choices of
the $u_{ii}$ is always sufficiently general for the BKK-bound to hold
with equality.

\begin{prb}
Is the ED degree of a torus $(S^1)^m \cong G \subseteq \GL(V)$ always
{\em equal} to the normalised volume of the convex hull of the character
set $X_V \subseteq \ZZ^m$ appearing in $V_\CC$?
\end{prb}

\subsection*{Other reductive groups preserving the form}

After this more or less satisfactory result for tori, it is tempting
to hope that the ED degree of any group $G$ preserving the bilinear
form should be expressible in terms of the highest weights appearing
in the complexification $V_\CC$ as a $G_\CC$-module. After all, $G_\CC$
is then a reductive group and much is known about its representations.
By Proposition~\ref{prop:degree} an upper bound to the ED degree is the
degree of a general orbit of $G_\CC$ in its action on $\End_\CC(V_\CC)$
by left multiplication. A formula for this degree is known by
\cite{Kazarnovskii87}; see also \cite[Theorem 8]{Derksen95}. However,
the space $\liea{g}^\bot$ is in general not sufficiently general for
that upper bound to be tight.

To test said hope, we have experimented with $G_\CC$ equal to the image
of $\SL_2(\CC)$ in its irreducible representation $V_\CC$ of highest
weight $m$ with $m$ even. Thus $V_\CC$ is the $m$-th symmetric power $S^m
\CC^2$ where $\CC^2$ is the standard representation of $G_\CC$. Since
$m$ is even, the bilinear form on $S^m \CC^2$ induced by the symplectic
$\SL_2(\CC)$-invariant form on $\CC^2$ is, indeed, an invariant symmetric
bilinear form. In this case, the formula in \cite[Theorem 8]{Derksen95}
evaluates to $m^3$. Below is a small table of ED degrees, in which we
could not yet find a pattern.

\begin{center}
\begin{tabular}{c|c|c|c|c|c}
$m$ & 0 & 2 & 4 & 6 & 8 \\
\hline
ED-degree of $\SL_2$ on $S^m \CC^2$ & 1& 4 & 40 & 156 & 400\\
\hline
$m^3$& 0 & 8 & 64 & 216 & 512
\end{tabular}
\end{center}

Note that the formula $m^3$ does not apply for $m=0$ since the
representation does not have a finite kernel.  Note also that the $4$
is consistent with Theorem~\ref{thm:Orthogonal} and the fact that the map
$\SL_2(\CC) \to \SL(S^2 \CC^2)$ has image $G_\CC$ equal to $\SO_3(\CC)$.

\begin{prb}
Determine the ED-degree of $\SL_2(\CC)$ on $S^m \CC^2$ with
$m$ even. More generally, find a formula for that ED degree
for any group $G_\CC$ on a representation $V_\CC$ with an
invariant symmetric bilinear form.
\end{prb}

\section{Groups not preserving the inner product}
\label{sec:notpreserving}

If $G$ does not preserve the inner product $(.|.)$, then it is much harder
to compute (or even estimate) the ED degree of $G$.  The following two
classes of groups illustrate this.

\subsection*{The special linear groups}
Consider the group 
\[ G:=\SL^{\pm}(V)=\{x \in \End(V) \mid \det x = \pm 1\}.\]
Its Lie algebra $\liea{g}$ is the
space of matrices with trace equal to zero. Given a real data matrix $u$,
the critical equations for the nearest $x\in G$ become
\[ \langle x^t(u-x), a \rangle =0 \quad \forall a\in \End(V),\tr(a)=0
\text{ subject to } \det x =1. \]
Since this equation must hold for all $a$ with $\tr(a)=0$ we see that
$x^t(u-x)$ must be of the form $cI$ for some $c\in \RR$, so that $u =
c x^{-t} + x$. From this expression for $u$ we find
\[ u^t u = (cx^{-1} + x^t)(cx^{-t} +x)= c^2 (x^t x)^{-1} + 2cI + x^t x. \]
Hence $s:=x^t x$ must be a symmetric matrix of determinant $1$ satisfying
\begin{equation} \label{eq:s}
u^t u = c^2 s^{-1} + 2cI + s = s^{-1} (c I + s)^2
\end{equation}
Conversely, if $s$ is a symmetric determinant-1 matrix satisfying this
equation, then we can set $x:= u^{-t}(cI + s)$. This matrix then satisfies
\[ x^t x = (cI+s) u^{-1}u^{-t} (cI+s) = (cI+s)^2 (u^tu)^{-1} = s, \]
where we have used that $s$ commutes with $u^tu$. This means that
$\det x = \pm 1$. We also find $x^t(u-x)=(cI+s)-s=cI$. Thus to compute
the ED degree of $G$ it suffices to count the symmetric matrices $s$
solving~\eqref{eq:s}.

Let $\mu_1,\ldots,\mu_n$ be the eigenvalues of $u^tu$. Since $u$ is
general, these are all distinct, and \eqref{eq:s} forces $s$ to be
simultaneously diagonalisable with $u^tu$. Thus we need only find the
eigenvalues $\lambda_1,\ldots,\lambda_n$ of $s$, where $\lambda_i$
corresponds to $\mu_i$. The equation~\eqref{eq:s} translates into
\[ \mu_i = c^2 \frac{1}{\lambda_i} +2c + \lambda_i, \qquad i=1,\ldots,n. \]
Multiplying by $\lambda_i$ and adding the condition $\lambda_1\cdots\lambda_n =1$, the system to solve becomes
\begin{equation} \label{eq:SLn}
\begin{cases}
f_i:=c^2 + (2c-\mu_i)\lambda_i + \lambda_i^2 &= 0, \qquad i=1,\ldots,n\\
\lambda_1\cdots\lambda_n &=1.
\end{cases}
\end{equation}
Substituting $\lambda_n:=(\lambda_1 \cdots \lambda_{n-1})^{-1}$ into
$f_n$, we find an ideal $I$ generated by $n$ equations $f_1,\ldots,f_n$
in the ring $\RR[\lambda_1^{\pm 1},\ldots,\lambda_{n-1}^{\pm
1},\mu_1,\ldots,\mu_n,c]$.  This ideal is prime, since the
equations can be read as defining the graph of a map from the
Cartesian product of an $(n-1)$-dimensional torus with coordinates
$\lambda_1,\ldots,\lambda_{n-1}$ with the affine line with coordinate $c$
to the affine space with coordinates $\mu_1,\ldots,\mu_n$. The ED degree
is the degree of this map. To determine it, we determine the intersection
$I \cap \RR[\mu_1,\ldots,\mu_n,c]$. For this, we eliminate
the $\lambda_i$ successively, as follows. For $i=0,\ldots,n-1$ define
$\lambda_{(i)}:=\lambda_1 \cdots \lambda_i$. Define
\[ R_n:=
c^2 \lambda_{(n-1)}^2
+(2c-\mu_n)\lambda_{(n-1)}
+1, 
\] 
which is just $f_n$ multiplied by $\lambda_{(n-1)}^2$.  Now recursively
define, for $i=1,\ldots,n-1$,
\[ R_{i}:=\Res_{\lambda_{i}}(R_{i+1},f_i), \]
where $\Res$ is the {\em resultant} given by the determinant of a suitable
Sylvester matrix. The first two are as follows:\\
\begin{align*} 
R_{n-1}&=
\det
\begin{bmatrix}
c^2 \lambda_{(n-2)}^2 & (2c-\mu_n)\lambda_{(n-2)} & 1 & 0\\
0 & c^2 \lambda_{(n-2)}^2 & (2c-\mu_n)\lambda_{(n-2)} & 1 \\
1 & (2c-\mu_{n-1}) & c^2 & 0\\
0 & 1 & (2c-\mu_{n-1}) & c^2\\
\end{bmatrix}\\
&= c^8 \lambda_{(n-2)}^4 + \cdots + 1,
\end{align*}
where the dots stand for terms of degrees strictly between
$0$ and $4$ in $\lambda_{(n-2)}$; and similarly
{\small
\begin{align*}
R_{n-2}&=
\det 
\begin{bmatrix}
c^8 \lambda_{(n-3)}^4 & . & . & . & 1 & 0 \\
0 & c^8 \lambda_{(n-3)}^4 & . & . & . & 1\\
1 & (2c-\mu_{n-2}) & c^2 & 0 & 0 & 0 \\
0 & 1 & (2c-\mu_{n-2}) & c^2 & 0 & 0 \\
0 & 0 & 1 & (2c-\mu_{n-2}) & c^2 & 0 \\
0 & 0 & 0 & 1 & (2c-\mu_{n-2}) & c^2 \\
\end{bmatrix}\\
&=c^{24} \lambda_{(n-3)}^8 + \cdots + 1.
\end{align*}
}
By induction, we find $R_i=c^{m_i} \lambda_{(i-1)}^{2^{n-i+1}} +
\cdots + 1$ where the remaining terms have $\lambda_{(i-1)}$-degree
strictly between zero and $2^{n-i+1}$. The exponents $m_i$ satisfy
the recursion
\[ m_i=2m_{i+1} + 2\cdot 2^{n-i} \]
and $m_n=2$. This is solved by $m_i=(n-i+1)2^{n-i+1}$. In particular,
we find that $m_1=n 2^n$. By induction one can prove that
$f_1,\ldots,f_{i-1},R_i$ generate the intersection
\[ I_i:=I \cap \RR[\lambda_1^{\pm
1},\ldots,\lambda_{i-1}^{\pm 1},\mu_1,\ldots,\mu_n,c] \]
and that, modulo $I_i$, the variable $\lambda_{i-1}$
can be expressed as a $\QQ$-rational function of
$\lambda_1,\ldots,\lambda_{i-2},\mu_1,\ldots,\mu_n,c$. This can be
used to show that for generic choices of the $\mu_i$ the degree-$n2^n$
equation $R_1$ in $c$ lifts to as many distinct solutions to the
system~\eqref{eq:SLn}. Thus we have proved the following theorem.

\begin{thm}
The ED degree of $\SL^{\pm}(V)$ equals $n2^n$, and the ED
degree of $\SL(V)$ equals $n2^{n-1}$. 
\end{thm}

The last statement follows from the fact that there exists an orthogonal
transformation of $\End(V)$ that takes the matrices with determinant $1$
into the matrices with determinant $-1$ and vice versa (e.g., in matrix
terms, multiplying the first column by $-1$).  Hence the two connected
components of $\SL^{\pm}(V)$ have the same ED degree.

The proof gives rise to the following algorithm for finding the closest
matrix in $G$ to a given real data matrix $u$: first diagonalise $u^t
u$ as
\[ u^t u = T \diag(\mu_1,\ldots,\mu_n) T^t, \]
where $T$ is a real orthogonal transformation and the $\mu_i$ are
positive. Then successively eliminate $\lambda_n,\ldots,\lambda_1$ as
above, using Sylvester matrices for the resultants $R_i$. Compute all
real roots $c$ of $R_1$. For each of these, compute the corresponding
$\lambda_1,\ldots,\lambda_n$ from the kernels of the Sylvester matrices:
since the data is sufficiently general, each of those kernels will
be one-dimensional and spanned by a vector of powers of the relevant
$\lambda_i$. Since all $\lambda_i$ are $\QQ$-rational functions of
$\mu_1,\ldots,\mu_n,c$, the $\lambda_i$ are, indeed, real.  Then construct
$s$ by
\[ s=T \diag(\lambda_1,\ldots,\lambda_n) T^t. \]
Finally, construct $x$ by 
\[ x =  u^{-t}(cI+s).\]
We have already verified that $x$ satisfies $x^t x=s$, so that the
$\lambda_i$ are necessarily positive, but this can also be seen directly
from \eqref{eq:SLn}.

It would be useful to know in advance which real root $c$ corresponds
to the closest matrix $x$. Experiments with the algorithm above suggests
that it may be the real root that is smallest in absolute value.

\begin{prb}
Is it true that the real root $c$ of $R_1$ of smallest absolute value
gives rise to the matrix $x \in \SL^\pm(V)$ that is closest to $u$?
\end{prb}

\section{The symplectic groups}

As a final case in our quest for ED degree of real algebraic groups we
fix an even $n=2m \in \NN$ and study the symplectic group
\[ \Sp_n:=\{x \in \RR^{n \times n} \mid x^t J x = J\}, \]
where $J$ has the block structure
\[ J=\begin{bmatrix} 0 & 1 \\ -1 & 0 \end{bmatrix}. \]
In the case of $\SL(V)$ the ED degree did not depend on the choice of an
inner product on $V$, because $\SL(V)$ acts transitively on inner products
(up to positive scalars). But for $\Sp_n$ the ED degree may well depend
on the relative position of the symplectic form given by $J$ and the
inner product. A general study of orbits of pairs of a symplectic and
a symmetric form is performed in \cite{Lancaster05,Dieudonne46}, based
on classical work by Kronecker. We choose the standard inner product.
This choice is rather special in the sense that the complexified
group $\Sp_n(\CC)$ intersects the complexified group $\O_n(\CC)$ in a
large group, containing a copy of the group $\GL_m(\CC)$. This is not
immediately clear from the chosen coordinates, but relative to the basis
\[
v_1:=\frac{e_1+ie_{m+1}}{\sqrt{2}},\ldots,v_m:=\frac{e_m+ie_{2m}}{\sqrt{2}}, 
v_{m+1}:=\frac{ie_1+e_{m+1}}{\sqrt{2}},\ldots,v_{2m}:=\frac{ie_m+e_{m}}{\sqrt{2}} \]
of $\CC^n$ the symplectic form still has Gram matrix $J$, while the
standard symmetric bilinear form on $\CC^n$ has Gram matrix
\begin{equation} \label{eq:Gram}
\begin{bmatrix} 0 & i \\ i & 0 \end{bmatrix}. 
\end{equation}
Now all complex matrices that relative to the basis of the $v_i$ have
the block structure
\[ \begin{bmatrix} g & 0 \\ 0 & g^{-T} \end{bmatrix} \]
lie both in $\Sp_n(\CC)$ and in $\O_n(\CC)$. This shows that we could
have chosen the symmetric form with Gram matrix (a scalar multiple of)
that in~\eqref{eq:Gram}, without changing the ED degree.

We have implemented the equations~\ref{eq:crit} and computed
the ED degree for very small values of $n$. The resulting table
is as follows:\\
\begin{center}
\begin{tabular}{l|c|c|c}
n & 2 & 4 & 6 \\
\hline
\text{ED degree of $\Sp_n$} & 4 & 24 & 544.\\
\end{tabular}\\
\end{center}
\ \\
The pattern might be that the answer is $2^{m^2}+2^{2m-1}$,
but we do not know how to prove this. 

\begin{prb}
Determine the ED degree of $\Sp_n$ for general even $n$.
\end{prb}

\bibliographystyle{alpha}
\bibliography{diffeq,draismapreprint}

\newcommand{\etalchar}[1]{$^{#1}$}
\def\cprime{$'$}
\begin{thebibliography}{DHO{\etalchar{+}}13}

\bibitem[Ber76]{Bernstein76}
D.N. Bernstein.
\newblock The number of roots of a system of equations.
\newblock {\em Funct. Anal. Appl. 9}, pages 183--185, 1976.

\bibitem[DHO{\etalchar{+}}13]{Draisma13c}
Jan Draisma, Emil Horobet, Giorgio Ottaviani, Bernd Sturmfels, and Rekha~R.
  Thomas.
\newblock The euclidean distance degree of an algebraic variety.
\newblock 2013.
\newblock Preprint available from \verb+http://arxiv.org/abs/1309.0049+.

\bibitem[Die46]{Dieudonne46}
Jean Dieudonn\'e.
\newblock Sur la r\'eduction canonique des couples de matrices.
\newblock {\em Bull. Soc. Math. France}, 74:130--146, 1946.

\bibitem[DK95]{Derksen95}
Harm {Derksen} and Hanspeter {Kraft}.
\newblock {Constructive invariant theory.}
\newblock In {\em {Alg\`ebre non commutative, groupes quantiques et invariants.
  Septi\`eme contact Franco-Belge, Reims, France, June 26--30, 1995}}, pages
  221--244. Paris: Soci\'et\'e Math\'ematique de France, 1995.

\bibitem[Hig89]{Higham89}
N.J. Higham.
\newblock Matrix nearness problems and applications.
\newblock In {\em Applications of matrix theory, Proc. Conf., Bradford/UK
  1988}, volume~22 of {\em Inst. Math. Appl. Conf. Ser., New. Ser.}, pages
  1--27, 1989.

\bibitem[Hor86]{Horn86}
B.K.P. Horn.
\newblock {\em Robot Vision}.
\newblock MIT Press and McGraw-Hill, 1986.

\bibitem[Kaz87]{Kazarnovskii87}
B.~Ya. Kazarnovski{\u\i}.
\newblock Newton polyhedra and {B}ezout's formula for matrix functions of
  finite-dimensional representations.
\newblock {\em Funktsional. Anal. i Prilozhen.}, 21(4):73--74, 1987.

\bibitem[Kel75]{Keller75}
Joseph~B. Keller.
\newblock Closest unitary, orthogonal and hermitian operators to a given
  operator.
\newblock {\em Math. Mag.}, 48:192--197, 1975.

\bibitem[LR05]{Lancaster05}
P.~Lancaster and L.~Rodman.
\newblock Canonical forms for symmetric/skew-symmetric real matrix pairs under
  strict equivalence and congruence.
\newblock {\em Linear Algebra Appl.}, 406:1--76, 2005.

\bibitem[OSS13]{Ottaviani13}
Giorgio Ottaviani, Pierre-Jean Spaenlehauer, and Bernd Sturmfels.
\newblock Algebraic methods for structured low-rank approximation.
\newblock 2013.
\newblock Preprint, available from \verb+http://arxiv.org/abs/1311.2376+.

\bibitem[Sch66]{Schonemann66}
P.H. Schonemann.
\newblock A generalized solution of the orthogonal procrustes problem.
\newblock {\em Psychometrika}, 31:1--10, 1966.

\end{thebibliography}

\end{document}